\documentclass{amsart}
\usepackage{amsmath,amssymb,amsthm}
\usepackage[pdftex]{graphicx}
\usepackage{graphicx}
\usepackage{longtable}
\usepackage{stmaryrd}
\usepackage{tikz}
\usetikzlibrary{calc}
\usetikzlibrary{decorations.pathmorphing}
\newtheorem{defi}{Definition}[section]
\newtheorem{theo}{Theorem}[section]
\newtheorem{prop}[theo]{Proposition}
\newtheorem{lemm}[theo]{Lemma}
\newtheorem{cor}[theo]{Corollary}
\newtheorem{rem}[theo]{Remark}
\newtheorem{notation}{Notation}[section]

\makeatletter
    
    \@addtoreset{equation}{section}
\makeatother
\makeatletter
  \newcommand{\subsubsubsection}{\@startsection{paragraph}{4}{\z@}%
    {1.0\Cvs \@plus.5\Cdp \@minus.2\Cdp}%
    {.1\Cvs \@plus.3\Cdp}%
    {\reset@font\sffamily\normalsize}
  }
  \makeatother
\allowdisplaybreaks[4]

\def\${|\!|\!|}
\newcommand{\rs}{\varodot}
\newcommand{\pl}{\varolessthan}
\newcommand{\pr}{\varogreaterthan}

\newcommand{\putsp}[2]{
\node (#1) at (#2) {};
}

\newcommand{\putdot}[4]{
\node (#1) [fill=white,circle] at (#2) {};
\fill (#1) circle (2pt);
\node (#1') [#3] at (#1) {$#4$};
}

\newcommand{\intputdot}[4]{
\node (#1) [fill=white,circle] at (#2) {};
\draw (#1) circle (2pt);
\node (#1') [#3] at (#1) {$#4$};
}

\newcommand{\putH}[4]{
\draw[very thick,->>] (#1)--(#2) node [#3,midway] {\tiny$#4$};
}

\newcommand{\puth}[4]{
\draw[->] (#1)--(#2) node [#3,midway] {\tiny$#4$};
}

\newcommand{\putv}[4]{
\draw[->,double] (#1)--(#2) node [#3,midway] {\tiny$#4$};
}

\newcommand{\putcv}[4]{
\draw[->,double] (#1)--(#2) node [#3,midway] {\tiny$#4$};
\path[fill=white, draw=black] ($(#1)!.5!(#2)$) circle (2pt);
}

\newcommand{\putdv}[4]{
\draw[->,densely dotted,double] (#1)--(#2) node [#3,midway] {\tiny$#4$};
}

\newcommand{\putcdv}[4]{
\draw[->,densely dotted,double] (#1)--(#2) node [#3,midway] {\tiny$#4$};
\path[fill=white, draw=black] ($(#1)!.5!(#2)$) circle (2pt);
}

\newcommand{\putddv}[4]{
\draw[->,densely dashed,double] (#1)--(#2) node [#3,midway] {\tiny$#4$};
}

\newcommand{\putcddv}[4]{
\draw[->,densely dashed,double] (#1)--(#2) node [#3,midway] {\tiny$#4$};
\path[fill=white, draw=black] ($(#1)!.5!(#2)$) circle (2pt);
}

\newcommand{\puthp}[5]{
\draw[->] (#1)--(#2) node [#3,pos=#4] {\tiny$#5$};
}

\newcommand{\putQ}[4]{
\draw[->>,very thick,decorate,decoration={snake,amplitude=0.5mm,segment length=1.5mm,post length=3mm}] (#1)--(#2) node [#3,midway] {\tiny$#4$};
}

\newcommand{\putq}[4]{
\draw[->,decorate,decoration={snake,amplitude=0.5mm,segment length=1.5mm,post length=1mm}] (#1)--(#2) node [#3,midway] {\tiny$#4$};
}

\newcommand{\putdq}[4]{
\draw[->,double,decorate,decoration={snake,amplitude=0.5mm,segment length=1.5mm,post length=1mm}] (#1)--(#2) node [#3,midway] {\tiny$#4$};
}

\newcommand{\putcdq}[4]{
\draw[->,double,decorate,decoration={snake,amplitude=0.5mm,segment length=1.5mm,post length=1mm}] (#1)--(#2) node [#3,midway] {\tiny$#4$};
\path[fill=white, draw=black] ($(#1)!.5!(#2)$) circle (2pt);
}

\newcommand{\putqp}[5]{
\draw[->,decorate,decoration={snake,amplitude=0.5mm,segment length=1.5mm,post length=1mm}] (#1)--(#2) node [#3,pos=#4] {\tiny$#5$};
}

\newcommand{\putch}[6]{
\draw[->] (#1)..controls(#2)and(#3)..(#4) node [#5,midway] {\tiny$#6$};
}

\newcommand{\putcH}[6]{
\draw[very thick,->>] (#1)..controls(#2)and(#3)..(#4) node [#5,midway] {\tiny$#6$};
}

\newsavebox{\boxI}
\sbox{\boxI}{\begin{tikzpicture}
\coordinate (A1) at (0,0);
\coordinate (A2) at ($(0,0.2)$);
\fill (A2) circle (1pt);
\draw (A1)--(A2);
\end{tikzpicture}}
\newcommand{\I}{{\usebox{\boxI}}}

\newsavebox{\boxY}
\sbox{\boxY}{\begin{tikzpicture}
\coordinate (A1) at (0,0);
\coordinate (A2) at ($0.6*(0.1,0.2)$);
\coordinate (A3) at ($0.6*(-0.1,0.2)$);
\coordinate (A4) at ($0.6*(0,-0.2)$);
\foreach \n in {2, 3} \fill (A\n) circle (0.9pt);
\draw (A3)--(A1)--(A2);
\draw (A1)--(A4);
\end{tikzpicture}}
\newcommand{\Y}{{\usebox{\boxY}}}

\newsavebox{\boxV}
\sbox{\boxV}{\begin{tikzpicture}
\coordinate (A1) at (0,0);
\coordinate (A2) at ($0.9*(0.1,0.2)$);
\coordinate (A3) at ($0.9*(-0.1,0.2)$);
\foreach \n in {2, 3} \fill (A\n) circle (1pt);
\draw (A3)--(A1)--(A2);
\end{tikzpicture}}
\newcommand{\V}{{\usebox{\boxV}}}

\newsavebox{\boxW}
\sbox{\boxW}{\begin{tikzpicture}
\coordinate (A1) at (0,0);
\coordinate (A2) at ($0.5*(-0.1,0.2)$);
\coordinate (A3) at ($0.5*(0.1,0.2)$);
\coordinate (A4) at ($0.5*(-0.2,0.4)$);
\coordinate (A5) at ($0.5*(0,0.4)$);
\coordinate (A6) at ($0.5*(0,-0.2)$);
\foreach \n in {3,4,5} \fill (A\n) circle (0.8pt);
\draw (A2)--(A1)--(A3);
\draw (A4)--(A2)--(A5);
\draw (A6)--(A1);
\end{tikzpicture}}
\newcommand{\W}{{\usebox{\boxW}}}

\newsavebox{\boxWc}
\sbox{\boxWc}{\begin{tikzpicture}
\coordinate (A1) at (0,0);
\coordinate (A2) at ($0.5*(-0.1,0.2)$);
\coordinate (A3) at ($0.5*(0.1,0.2)$);
\coordinate (A4) at ($0.5*(-0.2,0.4)$);
\coordinate (A5) at ($0.5*(0,0.4)$);
\foreach \n in {3,4,5} \fill (A\n) circle (0.9pt);
\draw (A2)--(A1)--(A3);
\draw (A4)--(A2)--(A5);
\draw (A3) to [out=90,in=0] (A5);
\end{tikzpicture}}
\newcommand{\Wc}{{\usebox{\boxWc}}}

\newsavebox{\boxIWc}
\sbox{\boxIWc}{\begin{tikzpicture}
\coordinate (A1) at (0,0);
\coordinate (A2) at ($0.5*(-0.1,0.2)$);
\coordinate (A3) at ($0.5*(0.1,0.2)$);
\coordinate (A4) at ($0.5*(-0.2,0.4)$);
\coordinate (A5) at ($0.5*(0,0.4)$);
\coordinate (A6) at ($0.5*(0,-0.2)$);
\foreach \n in {3,4,5} \fill (A\n) circle (0.9pt);
\draw (A2)--(A1)--(A3);
\draw (A4)--(A2)--(A5);
\draw (A6)--(A1);
\draw (A3) to [out=90,in=0] (A5);
\end{tikzpicture}}
\newcommand{\IWc}{{\usebox{\boxIWc}}}

\newsavebox{\boxWcc}
\sbox{\boxWcc}{\begin{tikzpicture}
\coordinate (A1) at (0,0);
\coordinate (A2) at ($0.5*(-0.1,0.2)$);
\coordinate (A3) at ($0.5*(0.1,0.2)$);
\coordinate (A4) at ($0.5*(-0.2,0.4)$);
\coordinate (A5) at ($0.5*(0,0.4)$);
\foreach \n in {3,4,5} \fill (A\n) circle (0.9pt);
\draw (A2)--(A1)--(A3);
\draw (A4)--(A2)--(A5);
\draw (A3) to [out=90,in=0] (A5);
\draw (A1) [fill=white] circle (1pt);
\end{tikzpicture}}
\newcommand{\Wcc}{{\usebox{\boxWcc}}}

\newsavebox{\boxK}
\sbox{\boxK}{\begin{tikzpicture}
\coordinate (A1) at (0,0);
\coordinate (A2) at ($1*(-0.1,0.1)$);
\coordinate (A3) at ($1*(0,0.2)$);
\fill (A3) circle (1pt);
\draw (A1)--(A2)--(A3);
\end{tikzpicture}}
\newcommand{\K}{{\usebox{\boxK}}}

\newsavebox{\boxB}
\sbox{\boxB}{\begin{tikzpicture}
\coordinate (A1) at (0,0);
\coordinate (A2) at ($0.6*(-0.2,0.2)$);
\coordinate (A3) at ($0.6*(0.2,0.2)$);
\coordinate (A4) at ($0.6*(-0.3,0.4)$);
\coordinate (A5) at ($0.6*(-0.1,0.4)$);
\coordinate (A6) at ($0.6*(0.1,0.4)$);
\coordinate (A7) at ($0.6*(0.3,0.4)$);
\foreach \n in {4,5,6,7} \fill (A\n) circle (0.9pt);
\draw (A2)--(A1)--(A3);
\draw (A4)--(A2)--(A5);
\draw (A6)--(A3)--(A7);
\end{tikzpicture}}
\newcommand{\Bo}{{\usebox{\boxB}}}

\newsavebox{\boxBc}
\sbox{\boxBc}{\begin{tikzpicture}
\coordinate (A1) at (0,0);
\coordinate (A2) at ($0.6*(-0.2,0.2)$);
\coordinate (A3) at ($0.6*(0.2,0.2)$);
\coordinate (A4) at ($0.6*(-0.3,0.4)$);
\coordinate (A5) at ($0.6*(-0.1,0.4)$);
\coordinate (A6) at ($0.6*(0.1,0.4)$);
\coordinate (A7) at ($0.6*(0.3,0.4)$);
\foreach \n in {4,5,6,7} \fill (A\n) circle (0.9pt);
\draw (A2)--(A1)--(A3);
\draw (A4)--(A2)--(A5);
\draw (A6)--(A3)--(A7);
\draw (A5) to [out=30,in=150] (A6);
\end{tikzpicture}}
\newcommand{\Bc}{{\usebox{\boxBc}}}

\newsavebox{\boxD}
\sbox{\boxD}{\begin{tikzpicture}
\coordinate (A1) at (0,0);
\coordinate (A2) at ($0.4*(-0.1,0.2)$);
\coordinate (A3) at ($0.4*(0.1,0.2)$);
\coordinate (A4) at ($0.4*(-0.2,0.4)$);
\coordinate (A5) at ($0.4*(0,0.4)$);
\coordinate (A6) at ($0.4*(-0.3,0.6)$);
\coordinate (A7) at ($0.4*(-0.1,0.6)$);
\foreach \n in {3,5,6,7} \fill (A\n) circle (0.8pt);
\draw (A2)--(A1)--(A3);
\draw (A4)--(A2)--(A5);
\draw (A6)--(A4)--(A7);
\draw (A1) [fill=white] circle (1pt);
\end{tikzpicture}}
\newcommand{\D}{{\usebox{\boxD}}}

\newsavebox{\boxC}
\sbox{\boxC}{\begin{tikzpicture}
\coordinate (A1) at (0,0);
\coordinate (A2) at ($1*(0.1,0.1)$);
\coordinate (A3) at ($1*(-0.1,0.1)$);
\coordinate (A4) at ($1*(0,0.2)$);
\foreach \n in {2, 4} \fill (A\n) circle (1pt);
\draw (A4)--(A3)--(A1)--(A2);
\draw (A1) [fill=white] circle (1pt);
\end{tikzpicture}}
\newcommand{\Co}{{\usebox{\boxC}}}

\newsavebox{\boxDa}
\sbox{\boxDa}{\begin{tikzpicture}
\coordinate (A1) at (0,0);
\coordinate (A2) at ($0.4*(-0.1,0.2)$);
\coordinate (A3) at ($0.4*(0.1,0.2)$);
\coordinate (A4) at ($0.4*(-0.2,0.4)$);
\coordinate (A5) at ($0.4*(0,0.4)$);
\coordinate (A6) at ($0.4*(-0.3,0.6)$);
\coordinate (A7) at ($0.4*(-0.1,0.6)$);
\foreach \n in {3,5,6,7} \fill (A\n) circle (0.8pt);
\draw (A2)--(A1)--(A3);
\draw (A4)--(A2)--(A5);
\draw (A6)--(A4)--(A7);
\draw (A1) [fill=white] circle (1pt);
\draw (A3) to [out=70,in=20] (A7);
\end{tikzpicture}}
\newcommand{\Da}{{\usebox{\boxDa}}}

\newsavebox{\boxDb}
\sbox{\boxDb}{\begin{tikzpicture}
\coordinate (A1) at (0,0);
\coordinate (A2) at ($0.4*(-0.1,0.2)$);
\coordinate (A3) at ($0.4*(0.1,0.2)$);
\coordinate (A4) at ($0.4*(-0.2,0.4)$);
\coordinate (A5) at ($0.4*(0,0.4)$);
\coordinate (A6) at ($0.4*(-0.3,0.6)$);
\coordinate (A7) at ($0.4*(-0.1,0.6)$);
\foreach \n in {3,5,6,7} \fill (A\n) circle (0.8pt);
\draw (A2)--(A1)--(A3);
\draw (A4)--(A2)--(A5);
\draw (A6)--(A4)--(A7);
\draw (A1) [fill=white] circle (1pt);
\draw (A3) to [out=90,in=0] (A5);
\end{tikzpicture}}
\newcommand{\Db}{{\usebox{\boxDb}}}

\newsavebox{\boxDc}
\sbox{\boxDc}{\begin{tikzpicture}
\coordinate (A1) at (0,0);
\coordinate (A2) at ($0.4*(-0.1,0.2)$);
\coordinate (A3) at ($0.4*(0.1,0.2)$);
\coordinate (A4) at ($0.4*(-0.2,0.4)$);
\coordinate (A5) at ($0.4*(0,0.4)$);
\coordinate (A6) at ($0.4*(-0.3,0.6)$);
\coordinate (A7) at ($0.4*(-0.1,0.6)$);
\foreach \n in {3,5,6,7} \fill (A\n) circle (0.8pt);
\draw (A2)--(A1)--(A3);
\draw (A4)--(A2)--(A5);
\draw (A6)--(A4)--(A7);
\draw (A1) [fill=white] circle (1pt);
\draw (A5) to [out=90,in=0] (A7);
\end{tikzpicture}}
\newcommand{\Dc}{{\usebox{\boxDc}}}

\def\a{\alpha}
\def\b{\beta}
\def\g{\gamma}
\def\de{\delta}
\def\e{\epsilon}

\def\k{\kappa}
\def\l{\lambda}
\def\s{\sigma}

\def\vp{\varphi}

\def\kpz{\cal{X}_{\text{\rm kpz}}}

\newcommand{\De}{\Delta}
\newcommand{\dl}{\partial}

\newcommand{\f}[2]{\frac{#1}{#2}}
\newcommand{\R}{\mathbb{R}}
\newcommand{\re}{\mathfrak{R}}
\newcommand{\T}{\mathbb{T}}
\newcommand{\w}[1]{\widetilde{#1}}
\newcommand{\Z}{\mathbb{Z}}

\newcommand{\dx}{\partial_x}
\newcommand{\dt}{\partial_t}

\newcommand{\bb}[1]{\mathbb{#1}}
\newcommand{\m}[1]{\mathbf{#1}}
\newcommand{\cal}[1]{\mathcal{#1}}
\newcommand{\fo}[1]{\widehat{#1}}

\newcommand{\sn}[3]{\|#1\|_{C_{#3}\cal{C}^{#2}}}
\newcommand{\hn}[4]{\|#1\|_{\cal{L}_{#4}^{#2,#3}}}

\newcommand{\hm}[5]{\|#1\|_{\cal{L}_{#5}^{#2,#3,#4}}}

\newcommand{\supp}{\mathop{\text{\rm supp}}}
\newcommand{\id}{\mathop{\text{\rm id}}}

\newcommand{\Ex}{\mathbb{E}}
\newcommand{\Px}{\mathbb{P}}
\def\paper{paper}

\begin{document}

\title[\resizebox{4.5in}{!}{Paracontrolled calculus and Funaki-Quastel approximation for the KPZ equation}]{Paracontrolled calculus and Funaki-Quastel approximation for the KPZ equation}
\author{Masato Hoshino}
\date{\today}
\address{The University of Tokyo, 3-8-1 Komaba, Meguro-ku, Tokyo 153-8914, Japan}
\email{hoshino@ms.u-tokyo.ac.jp}
\renewcommand{\include}[1]{}
\renewcommand\documentclass[2][]{}
\maketitle

\begin{abstract}
In this paper, we consider the approximating KPZ equation introduced by Funaki and Quastel \cite{FQ}, which is suitable for studying invariant measures. They showed that the stationary solution of the approximating equation converges to the Cole-Hopf solution of the KPZ equation with extra term $\f{1}{24}t$. On the other hand, Gubinelli and Perkowski \cite{GP} gave a pathwise meaning to the KPZ equation as an application of the paracontrolled calculus. We show that Funaki and Quastel's result is extended to nonstationary solutions by using the paracontrolled calculus.
\end{abstract}

\documentclass{amsart}

\section{Introduction}

The KPZ equation is the stochastic PDE
\begin{align}\label{2_intro:kpz eq}
\dt h(t,x)=\tfrac{1}{2}\dx^2h(t,x)+\tfrac{1}{2}(\dx h(t,x))^2+\dot{W}(t,x),\quad t>0,\ x\in\R,
\end{align}
where $\dot{W}$ is a space-time white noise, which is a centered Gaussian system with the covariance structure
\[\Ex[\dot{W}(t,x)\dot{W}(s,y)]=\de(t-s)\de(x-y).\]
We consider the equation \eqref{2_intro:kpz eq} on the torus $\T=\R/\Z$, equivalently on the interval $[0,1]$ with a periodic boundary condition.

The KPZ equation \eqref{2_intro:kpz eq} was introduced by Kardar, Parisi and Zhang \cite{KPZ} as a model for a growing interface represented by the height function $h$ with fluctuations. However, the equation \eqref{2_intro:kpz eq} is ill-posed. Indeed, we can expect that $h$ has a regularity $(\f{1}{2}-\delta)$ for every $\delta>0$ in the spatial variable, but this suggests that the non-linear term $(\dx h)^2$ would diverge. In order to cancel this singularity, we need to introduce the renormalized form of \eqref{2_intro:kpz eq}, which would be given by
\begin{align}\label{2_intro:reno kpz eq}
\dt h=\tfrac{1}{2}\dx^2h+\tfrac{1}{2}\{(\dx h)^2-\infty\}+\dot{W}.
\end{align}
By formally applying It\^o's formula, we can show that the solution $h$ of \eqref{2_intro:reno kpz eq} is given by the Cole-Hopf transform $h=\log Z$, where $Z$ is the solution of the stochastic heat equation with a multiplicative noise:
\begin{align}\label{2_intro:mSHE}
\dt Z=\tfrac{1}{2}\dx^2Z+Z\dot{W}.
\end{align}
We call $h_{\text{CH}}=\log Z$ the \emph{Cole-Hopf solution} of the KPZ equation.

In order to link the equation \eqref{2_intro:reno kpz eq} to the Cole-Hopf solution directly, we need to consider an approximation scheme. A natural approach is to replace $\dot{W}$ by a smeared noise $\dot{W}^\e(t,x)=(\dot{W}(t)*\eta^\e)(x)$ defined by a mollifier $\eta^\e=\e^{-1}\eta(\e^{-1}\cdot)$, where $\eta\in C_0^\infty(\R)$ is even and satisfies $\int\eta=1$, and to consider the equation
\begin{align}\label{2_intro:appro kpz eq}
\dt h^\e=\tfrac{1}{2}\dx^2h^\e+\tfrac{1}{2}\{(\dx h^\e)^2-C^\e\}+\dot{W}^\e,
\end{align}
where $C^\e=\int_\R\eta^\e(x)^2dx$. By applying It\^o's formula, we can show that $Z^\e=e^{h^\e}$ solves the equation
\[\dt Z^\e=\tfrac{1}{2}\dx^2Z^\e+Z^\e\dot{W}^\e.\]
It is easy to see that the solution $Z^\e$ converges to that of \eqref{2_intro:mSHE} as $\e\downarrow0$, therefore the solution $h^\e$ of \eqref{2_intro:appro kpz eq} also converges to the Cole-Hopf solution $h_{\text{CH}}$. For example, see Theorem 3.2 of \cite{BG}.

In order to study the invariant measures of the KPZ equation, \eqref{2_intro:appro kpz eq} is not a good approximation. Instead, Funaki and Quastel \cite{FQ} studied the approximation
\begin{align}\label{2_intro:fqkpz}
\dt \w{h}^\e=\tfrac{1}{2}\dx^2\w{h}^\e+\tfrac{1}{2}\{(\dx \w{h}^\e)^2-C^\e\}*\eta_2^\e+\dot{W}^\e,
\end{align}
where $\eta_2^\e=\eta^\e*\eta^\e$. They showed that the tilt process of the solution of \eqref{2_intro:fqkpz} has an invariant measure. Precisely, the distribution $\nu^\e$ of $\nabla(B*\eta^\e)=\{B*\eta^\e(x)-B*\eta^\e(y)\,;x,y\in\T\}$, where $B$ is a pinned Brownian motion on $\T$, is invariant under the tilt process $\nabla \w{h}^\e=\{\w{h}^\e(x)-\w{h}^\e(y)\,; x,y\in\T\}$. Furthermore, they also showed that the solution $\w{h}^\e$ with initial distribution $\nabla \w{h}^\e\sim\nu^\e$ converges to the process $h_{\text{CH}}(t,\cdot)+\f{1}{24}t$ in law sense.

Recently, Hairer \cite{Hairer13} gave a pathwise meaning to the KPZ equation based on the rough path theory. His method is a fixed point argument in a suitable Polish space independent of probability spaces and a control of several explicit stochastic processes. His work was extended to certain singular stochastic PDEs (e.g. dynamical $\Phi_d^4$ model, parabolic Anderson model, etc.), by Hairer's theory of \emph{regularity structures} \cite{Hai14}, and Gubinelli, Imkeller and Perkowski's \emph{paracontrolled calculus} \cite{GIP}.
In this \paper, we investigate the approximating equation \eqref{2_intro:fqkpz} by the paracontrolled calculus. As an application, we can show that the approximation \eqref{2_intro:fqkpz} works well for general initial values with positive regularity. Furthermore, the appearance of the constant $\f{1}{24}$ is easily explained as computed in Lemma 6.5 of \cite{Hairer13}.

Our main result is formulated as follows. We denote by $\mathcal{C}^\delta$ the Besov space $\mathcal{B}_{\infty,\infty}^\delta$ on $\T$, see Section \ref{2_section:para} for the precise definition.

\begin{theo}\label{2_main result}
Let $\vp\in C_0^\infty(\R)$ satisfy $\vp(0)=1$ and $\vp(x)=\vp(-x)$. Let $\eta=\cal{F}^{-1}\vp$ and consider the mollifier $\eta^\e=\e^{-1}\eta(\e^{-1}\cdot)$.
For every initial value $h_0\in\cal{C}^{0^+}=\cup_{\kappa>0}\cal{C}^\kappa$, there exists a survival time $T^\e\in(0,\infty]$ such that \eqref{2_intro:fqkpz} has a unique solution $\w{h}^\e$ on $[0,T]$ for every $T<T^\e$ and $\lim_{\e\downarrow0}T^\e=\infty$ in probability. Furthermore, $\w{h}^\e$ converges to the process $h(t)=h_{\text{\rm CH}}(t)+\frac1{24}t$ in $C((0,T],\cal{C}^{\f{1}{2}-\de})$ in probability for every $\de>0$ and $T<\infty$, where $h_{\text{\rm CH}}$ is the Cole-Hopf solution with initial value $h_0$.
\end{theo}

\begin{rem}
Precisely, the convergence $\w{h}^\epsilon\to h$ in probability considered here means that
$$
\Px(\|\w{h}^\epsilon-h\|_{C([t,T],\mathcal{C}^{\frac12-\delta})}>\lambda,\ T< T^\epsilon)+\Px(T\ge T^\epsilon)\to0
$$
for every $0<t<T$ and $\lambda>0$.
\end{rem}

This result is an extension of \cite{FQ} to non-stationary solutions and furthermore shows the convergence in probabilistically strong sense instead of law sense. This theorem follows from the following proposition. Its proof is given after the statement of Theorem \ref{2_drive1:convergence of X and wX}.

\begin{prop}\label{2_main result prop}
Let $h^\e$ and $\w{h}^\e$ be the solutions of the renormalized equations
\begin{align*}
&\dt h^\e=\tfrac{1}{2}\dx^2h^\e+\tfrac{1}{2}\{(\dx h^\e)^2-c^\e\}+\dot{W}^\e,\\
&\dt\w{h}^\e=\tfrac{1}{2}\dx^2\w{h}^\e+\tfrac{1}{2}\{(\dx\w{h}^\e)^2-\tilde{c}^\e\}*\eta_2^\e+\dot{W}^\e
\end{align*}
with common initial value $h_0\in\cal{C}^{0^+}$, where
\[c^\e=C^\e-\tfrac{1}{12}+\cal{O}(\e),\quad\tilde{c}^\e=C^\e.\]
Then both $h^\e$ and $\w{h}^\e$ converge to the process $h_{\text{\rm CH}}(t)+\frac1{24}t$.
\end{prop}

This \paper\ is organized as follows. In Section \ref{2_section:para}, we summarize some notations and results of paracontrolled calculus. In Section \ref{2_section:para kpz eq}, we give a pathwise meaning to the KPZ equation by paracontrolled calculus, and show the existence and uniqueness of the solution. In Section \ref{2_section:fqkpz}, we discuss the approximation \eqref{2_intro:fqkpz} by similar arguments to those in Section \ref{2_section:para kpz eq}. Section \ref{2_section:drive} is devoted to the probabilistic steps, i.e. the control of the driving stochastic processes.

\end{document}